\documentstyle[11pt]{amsart}

\newtheorem{prop}{Proposition}
\newtheorem{theo}{Theorem}
\newtheorem{Lemma}{Lemma}
\newtheorem{cor}{Corollary}

\newcommand{\na}{\nabla}
\newcommand{\om}{\omega}
\newcommand{\Om}{\Omega}
\newcommand{\la}{\lambda}
\newcommand{\La}{\Lambda}

\newcommand{\ka}{K{\"a}hler }

\newtheorem{conjecture*}{Conjecture}
\newtheorem{theorem*}{Theorem}
\newtheorem{question*}{Question}

\newcommand{\leftr}{[\hbox{\hspace{-0.15em}}[}
\newcommand{\rightr}{]\hbox{\hspace{-0.15em}}]}

\newcommand{\Lie}{\rm}

\newcommand{\llra}{\!\!\joinrel{\hbox to 30pt{\rightarrowfill}}}
\newcommand{\lllra}{\!\!\joinrel{\hbox to 50pt{\rightarrowfill}}}
\newcommand{\llllra}{\!\!\joinrel{\hbox to 60pt{\rightarrowfill}}}
\newcommand{\lllllra}{\!\!\joinrel{\hbox to 70pt{\rightarrowfill}}}
\newcommand{\llllllra}{\!\!\joinrel{\hbox to 75pt{\rightarrowfill}}}
\newcommand{\bx}{{\bf Q.E.D.}}

\title[rigidity of almost k{\"a}hler 4-manifolds]{Local rigidity of certain
classes
of almost k{\"a}hler 4-manifolds}
\author{VESTISLAV APOSTOLOV, JOHN ARMSTRONG AND TEDI DR\u{A}GHICI}
\thanks{The first author was supported in part by an FCAR, a PAFARC-UQAM,
and  by an NSERC grant. He is also member of EDGE, Research Training
Network HPRN-CT-2000-00101, supported by the European Human Potential
Programme. The first and the third authors were supported in part by
an NSF grant INT-9903302}
\address{Vestislav Apostolov \\ D\'epartement de math\'ematiques\\ UQAM
 \\ succursale Centre-ville  c.p. 8888 \\ Montr\'eal \\ H3C 3P8, Canada}
\email{apostolo@@math.uqam.ca}

\address{John Armstrong \\ 10, Alan Bullock Close \\ Oxford 0X4 1AU \\ UK}
\email{John.Armstrong@@madge.com}

\address{Tedi Dr\u{a}ghici \\ Department of Mathematics \\ Florida
International University \\
 Miami FL 33199 \\ USA}
\email{draghici@@fiu.edu}

\begin{document}

\begin{abstract}
We show that  any non-K{\"a}hler,
almost K{\"a}hler 4-manifold
for which both the Ricci
and the Weyl curvatures
have the same algebraic symmetries as they have for a K{\"a}hler
metric is locally isometric to the (only)
proper 3-symmetric 4-dimensional space \cite{gray,kow}.

\vspace{0.1cm}
\end{abstract}
\maketitle
\section{Introduction}
An {\it almost K{\"a}hler structure} on a manifold $M^{2n}$ is an
almost Hermitian structure
$(g, J, \Om)$ with a closed, and
therefore symplectic fundamental 2-form $\Om$.
If additionally the almost complex structure $J$
is integrable, then $(g, J, \Om)$ is a K{\"a}hler structure.
Almost K{\"a}hler metrics for which the almost complex structure is not
integrable
will be called {\it strictly} almost K{\"a}hler metrics.

Many efforts have been done in the direction of finding curvature
conditions on the metric which insure the integrability of the
almost complex structure. For example, an old, still open
conjecture of Goldberg \cite{Go} says that a compact almost
K{\"a}hler, Einstein manifold is necessarily K{\"a}hler.
Important progress was made by K. Sekigawa who proved that the
conjecture is true if the scalar curvature is non-negative
\cite{Se2}. The case of negative scalar curvature is still wide
open, despite of recent progress in dimension 4. Nurowski and
Przanowski \cite{NuP} and K.P.Tod \cite{Arm1,Tod} constructed
4-dimensional local examples of Einstein (in fact, Ricci flat),
strictly almost \ka manifolds. Thus, it is now known that
compactness must play an essential role, should the Goldberg
conjecture be true. In all these examples the structure of the
Weyl tensor is unexpectedly special --- the anti-self-dual part
of the Weyl tensor vanishes and the fundamental form is an
eigenform of the self-dual Weyl tensor (equivalently, $W^- =0$
and $W^+_2=0$, see below). Conversely, a recent result of
\cite{Arm1} states that any 4-dimensional strictly almost
K{\"a}hler, Einstein manifold is obtained by
Nurowski-Przanowski-Tod construction, provided that the
fundamental form is an eigenform of the Weyl tensor. It follows
that such a manifold can never be compact. Some other positive
partial results on the Goldberg conjecture in dimension 4 have
been obtained by imposing additional assumptions on the structure
of Weyl tensor, cf. \cite{AA,Arm0,Arm1,Arm,OS2}.

For an oriented four dimensional Riemannian manifold, it is well known
the ${\Lie{SO}}(4)$--decomposition of the Weyl tensor $W$ into its
self-dual and anti-self-dual parts, $W^+$ and
$W^-$. Moreover, for
every almost-Hermitian $4$-manifold $(M, g, J, \Om)$ the self-dual
part of the Weyl
tensor decomposes further under the action of the unitary group ${\Lie{U}}(2)$.
To see this, consider
$W^+$ as a trace-free, self-adjoint endomorphism of the bundle of self-dual
2-forms $\Lambda^+M$. Since $\Lambda^+M$ decomposes under
${\Lie{U}}(2)$ as ${\Bbb R}\Om \oplus \leftr \La^{0,2} M\rightr$, we
can write $W^+$ as a matrix with
respect to this block decomposition as follows:
\[
\left( \begin{array}{c|c}   \frac{\kappa}{6}       & W^+_2 \\ \hline
                          (W^+_2)^*    & W^+_3 - \frac{\kappa}{12}
{\rm Id}_{|\La^{0,2} M}
\end{array}
\right),
\]
where the smooth function $\kappa$ is the so-called conformal scalar curvature,
$W^+_2$ corresponds to the part of $W^+$ that interchanges the two
factors of the ${\Lie{U}}(2)$-splitting of $\La^+M$, and  $W^+_3$ is a
trace-free, self-adjoint endomorphism of the real vector
bundle $\leftr \La^{0,2}M  \rightr$ underlying the anti-canonical
bundle $\La^{0,2}M$. Also, the traceless part  of the
Ricci tensor ${\rm Ric}_0$ decomposes under $\Lie{U}(2)$ into two
irreducible components
 --- the invariant part  and
the anti-invariant part with respect to $J$, ${\rm Ric}_0^{\rm inv}$ and
${\rm Ric}_0^{\rm anti}$.
Correspondingly, there are several interesting types of
almost Hermitian $4$-manifolds, each imposing the vanishing of certain
${\rm U}(2)$-components of ${\rm Ric}_0$ and $W$, cf. \cite{TV}.

The curvature of a K{\"a}hler metric
$(g,J)$, for instance, satisfies any of the following three conditions:

\vspace{0.2cm}
\noindent
\begin{center}
(i) ${\rm Ric}_0^{\rm anti}=0$; \ (ii) $W^+_2=0$, and (iii) $W^+_3=0$.
\end{center}

\vspace{0.2cm}
\noindent
These three conditions are equivalent to the fact that
the curvature (considered as a ${\mathbb C}$-linear symmetric
endomorphism of the bundle of complex 2-forms) preserves the type
decomposition of 2-forms with respect to $J$, a property commonly
referred to as the {\it second Gray condition of the curvature}, cf.
\cite{Gr}.

Of course, the curvature of an arbitrary almost
K{\"a}hler metric may have none of these algebraic symmetries.
It is natural, therefore, to wonder if the integrability of the
almost complex structure is implied by the conditions (i)-(iii) above.
In \cite{ADK} and \cite{AD} an affirmative answer to this question is
given for  {\it compact} almost \ka 4-manifolds,
by using some powerful global arguments
coming from the Seiberg-Witten theory and Kodaira
classification of compact complex surfaces. One is then motivated to ask
what local rigidity, if any, do the conditions (i)-(iii) impose
on almost K\"ahler 4-manifolds. The goal of our paper is to answer
this question.

\vspace{0.2cm} We first provide a family of strictly almost \ka
4-manifolds satisfying, more generally, the conditions (i) and
(ii), see  Proposition \ref{prop1} below. Note that the strictly
almost K{\"a}hler, Ricci-flat flat examples of Nurowski,
Przanowski \cite{NuP} and Tod \cite{Arm1,Tod} satisfy (i) and
(ii) (but not (iii)), and our examples appear as a generalization
of Tod's construction \cite{Tod,Arm1}; instead of the
Gibbons-Hawking ansatz, we consider its generalized version
introduced by LeBrun in \cite{LeB}, and observe that appropriate
variable reductions lead to strictly almost \ka metrics with
$J$-invariant Ricci tensor and with special structure of the Weyl
tensor. While the Nurowski-Przanowski-Tod examples are just
particular metrics in this family, it turns out that for other
distinguished metrics the conditions (i)-(iii) are fulfilled.
Looking more carefully at the metrics satisfying conditions
(i)-(iii) from our family, one can further see that all of them
are, in fact, (locally) isometric to the unique 4-dimensional
proper (i.e. non-symmetric) {\it 3-symmetric space} described by
Kowalski \cite{kow} (see Section 4 below); as a homogeneous space
it is isomorphic to  $({\rm Isom}({\mathbb E}^2)\cdot
Sol_2)/SO(2)$ equipped with a left-invariant metric, or, by
introducing an invariant complex structure compatible with the
opposite orientation, it becomes isomorphic to the irreducible
homogeneous \ka surface  corresponding to the ${\bf
F_4}$-geometry of \cite{wall}. It might be also interesting to
note that this same example was discovered in yet a different
context by R. Bryant \cite{Br} (see also Remark 1).

Although one consequence of the existence of this example is that the
conditions
(i)-(iii) are not enough to insure the local integrability of an almost \ka
structure, we prove that, in fact, this is the only
such example in dimension four.

\begin{theo}\label{th1}
Any strictly almost \ka  4-manifold whose curvature satisfies
${\rm Ric}_0^{\rm anti}=0, \ W^+_2=0, \ W^+_3=0$ is locally isometric to
the (unique) 4-dimensional proper 3-symmetric space.
\end{theo}

\noindent {\it Remarks.} 1.--- It follows by Theorem \ref{th1}
and  the general theory of 3-symmetric spaces \cite{gray} that
any complete, simply connected strictly almost \ka 4-manifold
satisfying  the conditions (i)-(iii) is {\it globally} isometric
to the proper 3-symmetric 4-space.

2.--- Since
any 3-symmetric 4-space is almost \ka and satisfies (i)-(iii)
\cite{gray}, Theorem
\ref{th1} in turn provides a differential geometric proof of the
existence and the uniqueness of the proper 3-symmetric 4-space (see,
however,
\cite{kow} for more general results obtained by using Lie algebra
techniques).

3.---  Combining Theorem \ref{th1} with Wall's classification of
compact locally homogeneous complex surfaces \cite{wall}, one
sees that there are no {\it compact} strictly almost \ka
4-manifolds whose curvature satisfies the conditions (i)-(iii). This
provides an alternative proof of the integrability  result
in \cite{ADK} (see also Corollary 3 below).

\vspace{0.2cm}

Although our main goal of this paper is the study of almost \ka
4-manifolds which satisfy the three conditions (i)-(iii), Theorem
\ref{th1} is derived from the local classification of a larger
class of strictly almost K{\"a}hler 4-manifolds (Theorem
\ref{th2}), including as particular cases both the Einstein
metrics of \cite{NuP,Arm1} and the almost \ka 4-manifold
satisfying the conditions (i)-(iii) (see Remark 2). Our results
therefore generalize those in \cite{Arm1}.

The proof of our results relies on the strategy already developed
in \cite{Arm1} for finding out whether a given Riemannian metric
locally admits a compatible almost K{\"a}hler structure, which
allows us, as in \cite{Arm1}, to reduce the problem to an
integrable Frobenius system. However, the more general class of
almost \ka 4-manifolds that we consider in the current paper
leads to more involved proofs and makes the spinorial approach
invented in \cite{Arm1} somehow less adequate. We thus prefer to
use classical tensorial notations, which we hope will ease the
task of the reader in following the technical parts.

\vspace{0.2cm} The paper is organized as follows: In Sections 2
and 3, we prepare the necessary background of almost \ka geometry,
with a detailed analysis of the Riemannian curvature and its
covariant derivative, based on some representation theory. In
Section 4, we introduce our main examples of strictly almost \ka
4-manifolds satisfying conditions (i) and (ii), and describe
those which satisfy conditions (i)-(iii); we show that the latter
are isometric to the unique proper 3-symmetric 4-space. The last
section is devoted to the proof of our main result which is
stated in Theorem 2; Theorem 1 is then just a particular case.

\section{The curvature tensor of almost K{\"a}hler 4-manifolds}

Let $(M,g)$ be an oriented, 4-dimensional Riemannian manifold. The involutive
action of the Hodge
operator $*$ on the bundle of 2-forms $\Lambda ^2M$ induces the decomposition
$\Lambda^{2}M = \Lambda^{+}M \oplus  \Lambda^{-}M$ into the sub-bundles
of self-dual, resp. anti-self-dual 2-forms,
corresponding to the $+1$, resp. $-1$ eigenspaces of $*$.
We will implicitly identify vectors and covectors via the metric $g$ and,
accordingly,  a 2-form $\phi$ with the corresponding skew-symmetric
endomorphism of the tangent bundle $TM$, by putting: $g(\phi(X),Y) =
\phi(X,Y)$ for any vector fields $X,Y$.  Also, if $\phi,
\psi
\in TM^{\otimes 2}$, by $ \phi \circ \psi$ we understand the
endomorphism
of $TM$ obtained by the composition of the endomorphisms corresponding to
the two tensors. The inner product on $\La^2M$ induced by $g$ will be
denoted by $\langle \cdot , \cdot \rangle$, so as the induced norm
differs by a factor
$\frac{1}{2}$ from the usual tensor norm of $TM^{\otimes 2}$.

Considering the Riemannian curvature tensor $R$ as a symmetric
endomorphism of $\Lambda^2M$ we have the following  well known
$\Lie{SO}(4)$-splitting
\begin{equation}\label{so4}
R = \frac{s}{12}{\rm Id}_{| \La^2M} + \widetilde{{\rm Ric}_{0}} + W^{+} +
W^{-},
\end{equation}
where $s$ is the scalar curvature, $\widetilde{{\rm Ric}_{0}}$ is the the
Kulkarni-Nomizu
extension of the traceless Ricci tensor ${\rm Ric}_0$ to an endomorphism of
$\Lambda^2M$ (anti-commuting with $*$), and
$W^{\pm}$ are
respectively the self-dual and anti-self-dual parts of the Weyl tensor $W$.
The self-dual Weyl tensor  $W^{+}$
is viewed as a section of the bundle $S_{0}^2(\Lambda^{+}M)$ of  symmetric,
traceless  endomorphisms  of $\Lambda^{+}M$ (also considered as a sub-bundle of
the tensor product $\Lambda^{+}M \otimes \Lambda^{+}M$).

\vspace{0.2cm}

Let  $(M,g,J)$ be an almost Hermitian 4-manifold, {\it i.e.}, an
oriented Riemannian 4-manifold $(M,g)$ endowed with a
$g$-orthogonal almost complex structure $J$ which induces  the
chosen orientation of $M$. We denote by $\Om$ the corresponding
fundamental 2-form, defined by $\Om(X,Y) = g(JX,Y)$.
The action of $J$ extends to the cotangent bundle $\La^1M$ by putting
$(J\alpha)(X) = -\alpha(JX)$, so as to be compatible with the Riemannian
duality between $TM$ and $\La^1M$. This action defines  an {
involution}, $\imath_{J}$,
on $\Lambda ^2M$ by putting $\imath_J(\phi)(X,Y) = \phi(JX,JY)$, which
in turn gives rise to
the following orthogonal splitting of $\Lambda^{+}M$:
\begin{equation}\label{2}
\Lambda^{+}M = {\Bbb R} \Om \oplus \leftr \La^{0,2}M  \rightr  ,
\end{equation}
where  $\leftr \La^{0,2}M  \rightr $ denotes  the
bundle of $J$-anti-invariant real 2-forms, {\it i.e.}, the 2-forms
$\phi$ such that
$\imath_{J}(\phi)= - \phi$. Note that $\leftr \La^{0,2}M  \rightr $
is
the real underlying bundle of the anti-canonical bundle $(K_J)^{-1}=
\Lambda^{0,2}M$ of $(M, J)$; the induced complex structure $J$ on
$\leftr \La^{0,2}M  \rightr $ acts by $(J\phi)(X,Y)=-\phi(JX,Y)$.

Consequently, the vector bundle ${\cal W^+}=S_0^2(\Lambda^+M)$ of the
symmetric traceless
endomorphisms of $\Lambda ^+ M $ decomposes into the sum of three sub-bundles,
${\cal W}_{1}^+$, ${\cal W}_{2}^+$, ${\cal W}_{3}^+$, defined as follows,
see \cite{TV}:
\begin{enumerate}
\item[$\bullet$] ${\cal W}_1 ^+ = {\Bbb R} \times M $ is the sub-bundle of
elements
preserving the
decomposition (\ref{2}) and acting by homothety on the two factors;
hence it is
the trivial line
bundle generated by the element $ \frac{1}{8} \Om \otimes \Om - \frac{1}{12}
{\rm Id}_{| \La^+M}$.

\item[$\bullet$] ${\cal W}_2 ^+ = \leftr \La^{0,2}M  \rightr $ is the
sub-bundle of elements which
exchange the
two factors in (\ref{2}); the real isomorphism with $\leftr \La^{0,2}M
\rightr $ is seen
by identifying each
element $\phi$ of $\leftr \La^{0,2}M  \rightr $ with the element $
\frac{1}{2} (\Om \otimes
\phi + \phi \otimes \Om )$ of ${\cal W}_2 ^+$.

\item[$\bullet$] ${\cal W}_3 ^+ = S_0^2(\leftr \La^{0,2}M  \rightr )$ is
the sub-bundle of
elements preserving the splitting (\ref{2}) and acting trivially on
the first factor ${\Bbb R} \Om $.
\end{enumerate}
We then obtain the following $\Lie{U}(2)$-splitting of the Riemannian curvature
operator, cf. \cite{TV}:
\begin{equation}\label{u(2)}
R = \frac{s}{12} {\rm Id}_{| \La^2M} + ({\widetilde {{\rm Ric}_0}})^{\rm
inv} + ({\widetilde {{\rm Ric}_0}})^{\rm anti}
+
 W_1 ^+ + W_2 ^+ + W_3 ^+ + W^- ,
\end{equation}

\noindent
where $ ({\widetilde {{\rm Ric}_0}})^{\rm inv} $ and $ ({\widetilde {{\rm
Ric}_0}})^{\rm anti} $ are
the
Kulkarni-Nomizu extensions of the $J$-invariant and the
$J$-anti-invariants
parts of
the traceless Ricci tensor, respectively, and $W_i ^+$ are the projections
of $W^+$ on
the spaces ${\cal W}_i ^+, \; i = 1,2,3$. The component $W_1^+$
is given by
\begin{equation}\label{w^+1}
W_1^+ = \frac{\kappa}{8} \Om \otimes \Om - \frac{\kappa}{12} {\rm Id}_{|
\La^+M} ,
\end{equation}
where the smooth function $\kappa $ is the so called {\it conformal
scalar curvature} of $(g, J)$;
\begin{equation}\label{w^+2}
W^+_2=-\frac{1}{4}(\Psi\otimes \Omega + \Omega\otimes \Psi),
\end{equation}
for a section $\Psi$ of $\leftr \La^{0,2}M  \rightr $.\\
For any (local) section $\phi$ of $\leftr \La^{0,2}M
\rightr $ of square-norm 2,
the component in ${\cal W}^+_3$ is given by
\begin{equation}\label{w^+3}
W^+_3 = \frac{\lambda}{2}[\phi\otimes \phi - J\phi\otimes J\phi]
      + \frac{\mu}{2}[\phi\otimes J\phi + J\phi\otimes\phi],
\end{equation}
where  $\lambda $ and $\mu$ are (locally defined) smooth functions.

\vspace{0.2cm}

For any almost \ka structure $(g,J,\Om)$, the covariant derivative
$\na \Om$ of the fundamental form is identified with the {\it
Nijenhuis} tensor of $(M,J)$,  the obstruction for the
integrability of the almost complex  structure $J$. Moreover,
$\na \Om$ can be viewed as a section of the  real vector bundle
underlying $\La^{0,1}M\otimes \La^{0,2}M$, which allows us to write with
respect to any section $\phi$ of $\leftr \La^{0,2}M  \rightr$:
\begin{equation}\label{na-om}
\na\Om = a\otimes \phi - Ja \otimes J\phi.
\end{equation}
The 1-form $a$ satisfies $|\na \Om|^2 = 4|a|^2$, provided that
$\phi$ is of square-norm 2.  Consequently,
the covariant derivatives of $\phi$ and $J\phi$ are given
by
\begin{equation}\label{na-phi}
\na \phi = - a\otimes \Om + b\otimes J\phi; \ \na J\phi = Ja\otimes \Om -
b\otimes \phi,
\end{equation}
for some 1-form $b$.

Observe that
we have an $S^1$-freedom for
the choice of $\phi$ into the formulas (\ref{w^+3}) and (\ref{na-om}).
We shall refer to this as a {\it gauge
dependence} and any local section $\phi$ of $\leftr \La^{0,2}M
\rightr $ of square-norm 2
will be called a {\it gauge}.

\vspace{0.2cm}
\noindent
{\bf Convention.} From now on, $\phi$ will be assumed to be an
eigenform of $W^+_3$,
{\it i.e.}, the function $\mu$ in  (\ref{w^+3}) identically
vanishes.
\vspace{0.2cm}

Note that the above assumption can be locally arranged (for a
smooth gauge $\phi$ !) on the open dense subset of points, $x$,
where either
$W^+_3(x)
\neq 0$, or $W^+_3\equiv 0$ in the neighbourhood of $x$; however, by
continuity, all gauge independent properties will hold everywhere on
$M$.

\vspace{0.2cm}
Once the gauge $\phi$ is fixed as above, one can
determine the smooth functions
$\kappa$ and $\lambda$ and the 2-form $\Psi$
in terms of the 1-forms $a$ and $b$ and the 2-form
$\phi$, or, equivalently in terms of $2$-jets of $J$.
For that, we first make use of the {\it Weitzenb{\"o}ck formula} for
self-dual 2-forms, cf. {\it e.g.} \cite{bourg}:
\begin{equation}\label{weitz}
\Delta \psi = \na^*\na \psi + \frac{s}{3}\psi - 2W^+(\psi).
\end{equation}
Since the fundamental form $\Om$ is a self-dual, closed  2-form, it is
therefore  harmonic and (\ref{weitz}) implies
$$|\na \Om|^2 + \frac{2}{3}s-2\langle W^+(\Om), \Om \rangle = 0,$$
which, by (\ref{w^+1})--(\ref{w^+3}), is equivalent to
\begin{equation}\label{kappa}
\kappa - s = 6|a|^2 = \frac{3}{2}|\na \Om|^2.
\end{equation}

Formula (\ref{kappa})  shows that the smooth
function $\kappa - s$ is everywhere non-negative on $M$; it
vanishes exactly at the points where the Nijenhuis tensor is
zero. Observe also that applying (\ref{weitz}) to $\Om$ we involve
the 2-jets of $J$. Thus (\ref{kappa}) can be considered
as an ``obstruction''
to lifting the 1-jets of $J$ to 2-jets (see \cite{Arm1}), although
eventually it
takes form of a condition on the 1-jets.

\vspace{0.2cm}
\noindent
In order to express $W^+_2$ and $W^+_3$ we make use of the Ricci
identity
\begin{equation}\label{ricid}
(\na^2_{X,Y} - \na ^2_{Y,X})(\Om)(\cdot,\cdot) =
-R_{X,Y}(J\cdot,\cdot)-R_{X,Y}(\cdot,J\cdot).
\end{equation}
From (\ref{na-om}) we get
$$\na ^2|_{\Lambda^2 M}\Om =(da - Ja\wedge b)\otimes \phi -
(d(Ja) + a \wedge b)\otimes J\phi, $$
so, (\ref{ricid}) can be rewritten as
\begin{equation}\label{*}
da - Ja\wedge b = - R(J\phi); \ d(Ja) + a\wedge b = - R(\phi).
\end{equation}
Projecting on $\La^+M$ and using (\ref{u(2)})--(\ref{w^+3}) and (\ref{kappa}),
the equalities in (\ref{*})  give
\begin{eqnarray}\label{ricident}\label{lambda}
\lambda &=& - \frac{1}{2}\big{(}  |a|^2 - \langle da,J\phi \rangle
 + \phi(a,b)\big{)}; \\ \label{mu}
\mu &=& - \frac{1}{2}\big{(} \langle da,\phi\rangle   + J\phi(a, b)
 \big{)} =0;
\end{eqnarray}
\begin{equation}\label{Psi}
\Psi = \big{(} \langle d(Ja),\Om \rangle +
 \Om(a,b)\big{)}\phi +   \big{(} \langle da,\Om \rangle +
 g(a,b)\big{)} J\phi.
\end{equation}
We observe again that  the relations
(\ref{lambda})--(\ref{Psi}) are conditions on the
2-jets of the compatible almost \ka structure $J$, and can be viewed as a
further ``obstruction'' to lifting the $1$-jets to $2$-jets, see \cite{Arm1}.

Similarly, projecting formulae (\ref{*}) on $\La^-M$ we completely
determine the $J$-anti-invariant part of the Ricci tensor.
In order to determine its $J$-invariant part one
needs the 3-jets of $J$, involved  in
the Ricci identity for the Nijenhuis tensor (viewed as a section of
$\La^1M\otimes \La^2M$).
Writing the Ricci identity with respect
to $\na \Om$
is nothing but adding to (\ref{*}) one more relation coming from
$$(\na^2_{X,Y} - \na ^2_{Y,X})(\phi)(.,.) =
-R_{X,Y}(\phi.,.)-R_{X,Y}(.,\phi.).$$
Using (\ref{na-om}),(\ref{na-phi}) and
(\ref{u(2)})--(\ref{w^+3}) we eventually obtain
\begin{equation}\label{**}
db = a\wedge Ja -R(\Om) = a\wedge Ja - \frac{(s+2\kappa)}{12}\Om
-J\circ ({\rm Ric}_0^{\rm inv}) + \frac{1}{2}\Psi.
\end{equation}
The closed 2-form $db$ is gauge independent and is thus defined
on whole $M$; in fact, up to a factor $-\frac{1}{2\pi}$,  $db$ is a De
Rham representative of the first Chern class of $(M,J)$, see e.g.
\cite{H-F}.

Note that the relations (\ref{*}) and (\ref{**})
completely determine the Ricci tensor and the self-dual Weyl tensor of
$(M,g,J)$ in terms of the
$3$-jets of $J$. One can further see that
the remaining part of the curvature, the anti-self-dual Weyl tensor,
is determined by the $4$-jets of $J$. But we shall show in Section 5
that when the metric satisfies some additional properties, the relations
(\ref{*}) and (\ref{**}) are sufficient to write down
the whole Riemannian curvature of $g$. A careful
analysis of the above mentioned ``obstructions'' to lifting the $1$, $2$ and
$3$-jets of $J$ will eventually permit us to apply
the Frobenius theorem in
order to obtain the desired classification.

\section{Almost \ka 4-manifolds and Gray conditions. Preliminary results}

For a 4-dimensional almost Hermitian manifold, the relations
(i)--(iii) mentioned in the introduction are closely related to
the following conditions on the curvature defined by A. Gray
\cite{Gr} (not necessarily in the 4-dimensional context).
\newline $ (G_1) \; \; \; \; R_{XYZW} = R_{XYJZJW} $
 ;
\newline $ (G_2) \; \; \; \; R_{XYZW} - R_{JXJYZW} = R_{JXYJZW} + R_{JXYZJW}
$
 ;
\newline $ (G_3) \; \; \; \; R_{XYZW} = R_{JXJYJZJW}.$
\newline
Identity $(G_i)$ will be called the $i$-th Gray condition. Each imposes on
the curvature
of the almost Hermitian structure a certain degree of resemblance to that
of a \ka one.
A simple application of the first Bianchi identity
yields the implications $(G_1) \Rightarrow (G_2) \Rightarrow (G_3)$.
Also elementary is
the fact that a K{\"a}hler structure satisfies relation $(G_1)$ (hence, all
of the
relations $(G_i)$). Following \cite{Gr}, if ${\cal AK}$ is the class of
almost
\ka
manifolds, let ${\cal AK}_i$ be the subclass of manifolds whose curvature
satisfies
identity $(G_i)$. We have the obvious inclusions
$$ {\cal AK} \supseteq {\cal AK}_{3} \supseteq {\cal AK}_2 \supseteq
   {\cal AK}_1 \supseteq {\cal K} , $$
where ${\cal K}$ denotes the class of \ka manifolds. In \cite{Go}
it was observed that the equality $ {\cal AK}_1 = {\cal K}$ holds
   locally (this fact is an immediate consequence of
   (\ref{kappa})).

   From the examples of Davidov
and Mu\u{s}karov \cite{DM}, multiplied by compact \ka manifolds,
it follows that the inclusion $ {\cal AK}_2 \supset {\cal K}$ is
strict in dimension $2n \ge 6$, even in the compact case. This is
no longer true in dimension 4; it was proved in \cite{ADK} that
the equality ${\cal AK}_2 = {\cal K}$ holds for compact
4-manifolds (see also Corollary \ref{cor3} in Section 5 for a
partially different proof of this result).

\vspace{0.2cm}

Let us first observe that the conditions $(G_i)$ fit in with
the $\Lie{U}(2)$-decomposition (\ref{u(2)}) of the curvature in the
following manner:
\begin{Lemma}\label{lem1}
An almost Hermitian 4-manifold $(M, g, J)$ satisfies the property $(G_3)$
if and only if the Ricci tensor is $J$-invariant and $W_2 ^+ = 0$. It
satisfies $(G_2)$ if moreover $W_3 ^+ = 0$.
\end{Lemma}
\noindent
{\it Proof:} A consequence of (\ref{u(2)}), see \cite{TV}. $\bx$

\vspace{0.2cm}

Denote by ${\cal D}=\{ X \in T: \na_X \Om =0\}$ the {\it \ka nullity}
of $(g,J)$ and by
${\cal D}^{\perp}$ its $g$-orthogonal
complement. According to (\ref{na-om}), ${\cal D}$ is
$J$-invariant at every point and has rank 4 or 2,
depending on whether or not the Nijenhuis tensor $N$ vanishes at that point.
As an easy consequence of (\ref{*}), we have the following useful
observation:
\begin{Lemma}\label{lem2} A non-K{\"a}hler, almost \ka 4-manifold with
  $J$-invariant
Ricci tensor belongs to the class ${\cal AK}_3$ if and only if the \ka nullity
${\cal D}$  is a rank 2 involutive distribution on the open set of points
where the Nijenhuis tensor does not vanish.
\end{Lemma}
\noindent
{\it Proof:}
Let $\{ B,JB\}$ be any (local)
orthonormal frame of ${\cal D}$ and let
$\{A, JA\}$ be an orthonormal frame of ${\cal D}^{\perp}$, so that $A$
and $JA$ are the dual orthonormal frame of $\{ a,Ja\}$, see
(\ref{na-om}). Then the fundamental form can be written as
\begin{equation}\label{om}
\Om = A\wedge JA + B\wedge JB.
\end{equation}
By (\ref{*}) we see that ${\cal D}$
is involutive if and only if
\begin{equation}\label{temp1}
R(\phi)(B,JB)=0, \ \ \ R(J\phi)(B,JB)=0.
\end{equation}
On the other hand, as the Ricci tensor is $J$-invariant, it follows by
(\ref{u(2)})--(\ref{w^+3}) and (\ref{om}):
$$R(\phi)(B,JB)=-\frac{1}{4}\langle \Psi, \phi \rangle; \
R(J\phi)(B,JB)=-\frac{1}{4}\langle \Psi, J\phi \rangle,$$
{\it i.e.},  according to (\ref{temp1}), we obtain that ${\cal D}$ is
involutive if and only if $W^+_2=0$ (see (\ref{w^+2})). The claim now
follows by
Lemma \ref{lem1}. $\bx$

\vspace{0.2cm}

We shall further use the following refined version
of the differential Bianchi identity \cite{AA}:

\begin{Lemma}\label{lem3}{\rm {\bf (Differential Bianchi
      identity)}}
Let $(M, g, J)$ be an almost K{\"a}hler 4-manifold
in the class ${\cal AK}_3$. Then the following relations hold:
\begin{equation}\label{lem3-1}
d(\kappa -s) =  - 12\la J\phi(a);
\end{equation}
\begin{equation}\label{lem3-2}
{\rm Ric}_0(a)= \frac{\kappa}{4}a  + 2\la\phi(b) - J\phi(d\la);
\end{equation}
\begin{equation}\label{lem3-3}
\Delta (\kappa - s) = -\frac{\kappa}{2}(\kappa - s) - 24\la^2 + 12
{\rm Ric}_0(a,a).
\end{equation}
\end{Lemma}
\noindent
{\it Proof:} The co-differential $\delta W^+$ of the self-dual Weyl
tensor of $(M, g)$
is a section of the rank 8 vector bundle
$ {\cal V} = {\rm Ker} ({\it trace} : \La^1M \otimes \La^+M  \rightarrow
\La^1M ),$ where the trace  is defined by $ { trace} (\alpha \otimes
\phi) =
\phi(\alpha) $ on decomposed elements. For every almost-Hermitian
4-manifold the vector bundle
${\cal V}$ splits as
${\cal V} = {\cal V}^+ \oplus {\cal V}^-$, see \cite{AG}, where
${\cal V}^+$
is identified with the cotangent bundle $\La^1 M$ by
\begin{equation}
\label{alpha}
 \La^1M \ni \alpha \mapsto A = J\alpha \otimes \Om -
\frac{1}{2}\sum_{i=1}^{4}e_{i} \otimes (\alpha \wedge e_{i} - J\alpha \wedge
Je_{i}),
\end{equation}
while $\cal V^{-}$ is identified (as a real vector bundle)  with
$\Lambda^{0,1}M \otimes \La^{0,2}M$. For any gauge
$\phi$ the vector bundle ${\cal V}^-$
can be again identified with $\La^1 M$  by
\begin{equation}\label{beta}
\La^1M \ni \beta \mapsto B = J\beta \otimes \phi + \beta\otimes
J\phi.
\end{equation}
We denote by  $(\delta W^{+})^{+}$, resp. $(\delta W^{+})^{-}$, the
component of $\delta W^{+}$ on  $\cal V^{+}$, resp. on  $\cal
V^{-}$, and,  for any gauge $\phi$
satisfying the Convention of Section 2 we consider the corresponding 1-forms
$\alpha$ and $\beta$. By (\ref{alpha}),
(\ref{beta}) and (\ref{w^+1})--(\ref{w^+3}) one directly calculates:
\begin{equation}\label{alpha1}
\alpha= -\frac{1}{2}J\langle \delta W^+, \Om \rangle=
-\frac{d\kappa}{12}  -
\la J\phi(a);
\end{equation}
\begin{eqnarray}\label{beta1}
\beta &=& \frac{1}{2}\big{(}-J\langle \delta W^+, \phi \rangle +
      \frac{1}{2}\phi\langle \delta W^+, \Om \rangle \big{)}\\ \nonumber
      &=& -\frac{\kappa}{8}a +\la\phi(b) - \frac{1}{2}J\phi(d\la).
\end{eqnarray}
Recall that the {\it Cotton-York tensor} $C$ of $(M, g)$ is  defined by:
$$ C_{X,Y,Z} = \frac{1}{2} \Big[ \na_{Z} (\frac{s}{12} g + {\rm Ric}_0) (Y,X) -
\na_{Y} (\frac{s}{12} g + {\rm Ric}_0)(Z,X) \Big],$$
for any vector fields $X,Y,Z$. Considering $C$ as a 2-form with values
in $\La^1M$, the {\it second Bianchi identity} reads as
$\delta W = C$. In dimension 4 we have also the
``half'' Bianchi identity
\begin{equation}\label{half}
\delta W^+ = C^+,
\end{equation}
\noindent
where $C^+$ denotes the self-dual part of $C_X$, $X \in TM$.
When the Ricci tensor is  $J$-invariant, we
make use of (\ref{half}) to give an equivalent expression for the
1-forms $\alpha$ and $\beta$ in terms of the Ricci tensor and the
1-form $a$. According to (\ref{alpha}) we get
$$ \alpha(X) = -\frac{1}{2}J\langle C^+,\Om \rangle = - \frac{1}{4}
\sum_{i=1}^{4} \na_{e_{i}} (\frac{s}{12} g + {\rm Ric}_0)(Je_i, JX) = $$
$$ = - \frac{1}{4} \Big[ \frac{ds}{12}(X) - (\delta {\rm Ric}_0)(X) +
\sum_{i=1}^{4} {\rm Ric}_0(e_i, J(\na_{e_i} J)(X))  \Big] =$$
$$ = - \frac{1}{4} \Big[ \frac{ds}{3}(X) +
\sum_{i=1}^{4} {\rm Ric}_0(e_i, J(\na_{e_i} J)(X))\Big].$$
Using (\ref{na-om})  and the fact that the Ricci tensor is
$J$-invariant,
we obtain
$$ \sum_{i=1}^{4} {\rm Ric}_0(e_i, J(\na_{e_i} J)(X)) = 0,$$
and then
\begin{equation}\label{bianchi-a}
\alpha = -\frac{ds}{12}.
\end{equation}
Regarding the component of $C^+$ in ${\cal V}^-$, we have by (\ref{beta}):
$$\beta = \frac{1}{2}\big{(}-J\langle C^+, \phi \rangle
+\frac{1}{2}\phi\langle C^+,\Om \rangle \big{)}.$$
To compute $J\langle C^+, \phi \rangle $ we proceed in the same way as
computing $J\langle C^+, \Om \rangle$; instead of $J$ we
consider  the almost complex structure $I_{\phi}$
whose K{\"a}hler form is $\phi$. Observe that ${\rm Ric}_0$ is now
$I_{\phi}$-anti-invariant. By (\ref{na-om}),
(\ref{na-phi}) and (\ref{bianchi-a}) we eventually get
\begin{equation}\label{bianchi-b}
\beta = -\frac{1}{2}{\rm Ric}_0(a).
\end{equation}
Comparing (\ref{bianchi-a}) and (\ref{bianchi-b}) with (\ref{alpha1})
and (\ref{beta1}) we obtain the equalities (\ref{lem3-1}) and
(\ref{lem3-2}). Finally, taking co-differential of
 both sides of (\ref{lem3-1})
and using (\ref{lem3-2}) and (\ref{kappa}) we derive
\begin{eqnarray}\nonumber
\Delta (\kappa - s) &=& -12 J\phi(d\la, a) - 12\la \delta(J\phi(a)) \\
\nonumber
                    &=& 12{\rm Ric}_0(a,a) -\frac{\kappa}{2}(\kappa -s) + \\
\nonumber
                    & & 12\la \big(2\phi(a,b) - \langle d a, J\phi
\rangle + \delta (J\phi)(a)\big).
\end{eqnarray}
By (\ref{lambda}) and (\ref{na-phi}) we calculate
$$ 12\la \big(2\phi(a,b) - \langle da, J\phi \rangle + \delta
(J\phi)(a)\big) = -24\la^2,$$
and we reach the equality (\ref{lem3-3}). $\bx$

\vspace{0.2cm}
\noindent
We have the following consequence of Lemma \ref{lem3}
(see also \cite[Prop.2]{AD} and \cite[Prop.4]{Jel}):
\begin{cor}\label{cor1} A 4-dimensional almost \ka structure $(g, J,
\Om)$ in the
class ${\cal AK}_3$ belongs to ${\cal AK}_2$
if and only if the norm of $\na \Om$ is constant.
Moreover, if
$(g,J, \Om)$ is an ${\cal AK}_2$, non-\ka structure,
then the traceless Ricci tensor ${\rm Ric}_0$ is given by
$$ {\rm Ric}_0=\frac{\kappa}{4}[-g^{\cal D} + g^{{{\cal D}}^{\perp}}],$$
where $g^{\cal D}$ (resp. $g^{{\cal D}^{\perp}}$) denotes the restriction
of $g$ on  ${\cal D}$ (resp. on ${\cal D}^{\perp}$).
\end{cor}
\noindent
{\it Proof:}
According to (\ref{kappa}), we have
$|\na \Om|^2=\frac{\kappa - s}{6}$. We then get
by Lemma \ref{lem3} the equality $d(|\na \Om|^2) = -2\lambda J\phi(a),$
and the first part of the claim
follows by Lemma \ref{lem1} and (\ref{w^+3}). Since $W^+_3 \equiv 0$
({\it i.e.} $\la \equiv 0$ according  to (\ref{w^+3})),
the second relation stated in Lemma \ref{lem3} reads as
${\rm Ric}_0(a)=\frac{\kappa}{4}a.$
As ${\rm Ric}_0$ is symmetric traceless and $J$-invariant tensor,  in
the case when
$(g,J)$ is not K{\"a}hler the
expression above implies  the second part of the corollary. $\bx$

\section{Examples of almost \ka 4-manifolds satisfying Gray conditions}

\subsection{3-symmetric spaces} In this subsection we briefly describe
an already known example of strictly almost K\"ahler 4-manifold satisfying
the condition ($G_2$). This example comes from works of Gray
\cite{gray} and Kowalski \cite{kow} on {\it 3-symmetric spaces} and we
refer to their papers for more details on the subject.

A Riemannian 3-symmetric space is a manifold $(M,g)$ such that for each point
$p \in M$ there
exists an isometry $\theta_p  : M \rightarrow M$ of order 3
(i.e. $\theta_p^3 = 1$), with $p$ as an isolated fixed point.
Any such manifold has a naturally defined (canonical) $g$-orthogonal
almost complex structure $J$, and we further require that each $\theta_p$
is a holomorphic map
with respect to $J$. Moreover, the canonical almost Hermitian structure
$(g,J)$ of a 3-symmetric space always satisfies the second Gray condition
and, in dimension 4, is automatically almost K\"ahler
(it is K\"ahler if and only if the manifold is Hermitian symmetric,
see \cite{gray}).
It only remains the question of whether there exists a 4-dimensional
example of a 3-symmetric space with a non-integrable almost complex
structure (we shall call this a {\it proper} 3-symmetric space).
This is solved by Kowalski, who constructs such an example
and, moreover, shows that this is the only proper 3-symmetric space
in dimension 4 (in fact, this is the only proper {\it generalized}
symmetric space in dimension 4, \cite[Theorem VI.3]{kow}).
Explicitly, up to a homothety, Kowalski's example
is defined on ${\Bbb R}^4=\{(u_1,v_2,u_2,v_2 )\}$ with the metric
\begin{eqnarray} \label{kowalm}
g &=& \Big(-u_1 + \sqrt{u_1^2 + v_1^2 + 1}\Big)du_2^2 +
\Big(u_1 + \sqrt{u_1^2 + v_1^2 + 1}\Big)dv_2^2 \\
\nonumber
  & & - 2v_1 du_2\odot dv_2 \\ \nonumber
  & & + \frac{1}{(u_1^2 + v_1^2 +1)}\Big[(1+ v_1^2)du_1^2 +
(1 + u_1^2)dv_1^2
   - 2u_1v_1~ du_1 \odot dv_1 \Big],
\end{eqnarray}
where as usually $\odot$ stands for symmetric tensor products.

\subsection{Generalized Gibbons-Hawking Ansatz}

We now present a different and more general approach
of obtaining examples of almost \ka 4-manifolds satisfying Gray conditions
$(G_3)$ and $(G_2)$, which is based on the idea of generalizing Tod's
construction
of Ricci-flat strictly almost  \ka
4-manifolds \cite{Arm1,Tod}. For this purpose, we consider
instead of the Gibbons-Hawking ansatz, its
generalized version, introduced by LeBrun \cite{LeB} to construct
scalar-flat \ka surfaces. Following \cite{LeB}, let $w>0$ and $u$ be
smooth real-valued
functions on an open, simply-connected set $V \subset {\Bbb R}^3=\{(x,y,z)\}$,
which satisfy
\begin{equation}\label{1}
w_{xx} + w_{yy} + (we^u)_{zz}=0.
\end{equation}
Let $M = {\Bbb R}\times V$ and $\om$ be a 1-form on $M$
non-vanishing when restricted to the ${\Bbb R}$-factor and
determined (up to gauge equivalence) by
\begin{equation}\label{dom}
d\om = w_xdy\wedge dz + w_y dz\wedge dx + (we^u)_z
dx\wedge dy .
\end{equation}
It is shown in \cite{LeB} that the metric
\begin{equation}\label{g}
g = e^{u}w(dx^2 + dy^2) + w dz^2 + w^{-1}\om^2
\end{equation}
admits a \ka structure $I$, defined by its fundamental form
\begin{equation}\label{I}
\Om_{I}=dz\wedge \om + e^{u}w dx \wedge dy.
\end{equation}
Moreover, if we denote by $\frac{\partial}{\partial t}$
the dual vector field of $w^{-1}\om$ with respect to $g$,
then $\frac{\partial}{\partial t}$ is
Killing and preserves $I$.  Conversely, every \ka metric admitting a
hamiltonian Killing field locally arises by this construction \cite{LeB}.

\vspace{0.2cm}
\noindent
Besides the \ka
structure $I$, we shall consider the almost Hermitian structure $J$
whose fundamental form is
\begin{equation}\label{J}
\Om_J= - dz\wedge \om + e^{u}w dx \wedge dy.
\end{equation}
Clearly, the almost complex structures $I$ and $J$ commute
and yield different orientations on $M$.
Our  objective  is the following generalization of \cite{Tod}:
\begin{prop}\label{prop1} Let $w>0$ and $u$ be smooth functions
satisfying (\ref{1}). Then
the almost Hermitian structure $(g,J)$
defined via (\ref{g}) and (\ref{J}) is almost \ka if and only if
$u$ and $w$ satisfy
\begin{equation}\label{2'}
(e^{u}w)_z=0.
\end{equation}
It is  \ka if moreover $w$ does not depend on $x$ and $y$.

Furthermore, the following are true:
\begin{enumerate}
\item[{\rm (i)}] The almost
Hermitian manifold $(M,g,J)$ is non-\ka and belongs to ${\cal AK}_3$
if and only if $w$ is a non-constant, positive  harmonic
function of $x$ and $y$, and $u(x,y)$ is any function defined on
$U = V \cap {\Bbb R}^2$.
\item[{\rm (ii)}] The manifold $(M,g,J)$  belongs to ${\cal
AK}_2$ if and only if, in addition, $w$ has no critical values on $U$ and
$u$ is
given by
\end{enumerate}
\begin{equation}\label{u}
u= \ln(w_x^2 + w_y^2) - 3 \ln w + const.
\end{equation}
\end{prop}

\noindent
{\bf Remark 1.} (a) If $w$ is a non-constant harmonic function of $(x,y)$,
then the holomorphic function $h$ of $x+iy$ such that ${\rm Re}(h)=w$
can be used as a holomorphic coordinate in place of $x+iy$. Up to a change
of the smooth function $u$ and the transversal coordinate $t$, the
metrics described in Proposition \ref{prop1}(i) are then all
isometric to
\begin{equation}\label{ak3}
g = e^ux(dx^2 + dy^2) + xdz^2 + \frac{1}{x}(dt + ydz)^2,
\end{equation}
and therefore is defined on $M =\{(x,y,z,t) \in {\Bbb R}^4, x>0\}$
for any smooth function $u$ of $(x,y)$.
It is easily checked \cite{LeB} that the Ricci tensor
of the metrics (\ref{ak3}) has two vanishing eigenvalues while the
scalar curvature $s$ is given by
$s=\frac{u_{xx} + u_{yy}}{4xe^u}.$
It thus follows that the Ricci-flat Tod's examples are obtained precisely
when $u$ is a harmonic function.

(b) Concerning the metrics given in Proposition \ref{prop1}(ii),  by
(\ref{u})  we obtain
$e^u = const.\frac{1}{x^3}$, so that
(up to homothety of $(z,t)$) all these metrics
are  homothetic  to
\begin{equation}\label{canonic}
g = \frac{dx^2}{x^2} + \frac{1}{x^2}\sigma_1^2 + x\sigma_2^2 +
\frac{1}{x}\sigma_3^2,
\end{equation}
where $\sigma_1=dy~; \sigma_2=dz~; \sigma_3=dt + ydz$ are the standard
generators of the Lie algebra of the three dimensional Heisenberg group
${\rm Nil}^3$. It turns out that (\ref{canonic}) defines  a complete
metric, in fact,
a homogeneous one which is nothing else than the (unique)
proper 3-symmetric metric (\ref{kowalm}) mentioned in Sec. 4.1.
To see this directly, one should do the change of variables
\begin{equation}\label{chvar}
 u_1 = \frac{x^2 + y^2 - 1}{2x} , \; \; v_1 = - \frac{y}{x} , \; \; u_2 = t
, \; \; v_2 = z ,
\end{equation}
and after a straightforward calculation it can be seen that
the metric of Kowalski defined by (\ref{kowalm}) reduces exactly to
(\ref{canonic}).
In fact, we were motivated to look for and were able to find this change of
variables
only {after} we realized that one must have
the uniqueness stated in Theorem 1 (see also Remark 4).

(c) One can easily write down the whole Riemannian curvature of
the metric (\ref{canonic}):  it turns out that it is completely
determined by the (constant) scalar curvature $s=\frac{u_{xx} +
u_{yy}}{4xe^u}=-\frac{3}{4}$. Indeed, it is easily checked that
the conformal scalar curvature (which determines $W^+$)  is equal
to $\frac{3}{4}$, the Ricci tensor has constant eigenvalues
$(0,0,-\frac{3}{8},-\frac{3}{8})$, and as $g$ is K{\"a}hler with
respect to $I$ (see (\ref{I})), the anti-self-dual Weyl tensor is
also determined by $s$, see {\it e.g.} \cite{Ga}. The  metric
(\ref{canonic}) with its negative K{\"a}hler structure $I$
provide therefore a  non-symmetric, homogeneous K{\"a}hler
surface which corresponds to the ${\bf F_4}$-geometry of
\cite{wall}; it is thus a complete irreducible \ka metric with
two distinct {\it constant} eigenvalues of the Ricci tensor. From
this point of view, the metric (\ref{canonic}) has been
independently discovered by R. Bryant in \cite{Br}. Remark that
many others (non-homogeneous in general) K\"ahler metrics of
constant eigenvalues of the Ricci tensor arise from (\ref{ak3}),
provided that $u$ is a smooth solution to the elliptic equation
$$u_{xx} + u_{yy} = 4sxe^u,$$
where $s$ is a non-zero constant, the scalar curvature of the metric.

\vspace{0.1cm}
\noindent
{\it Proof of Proposition \ref{prop1}:}
By (\ref{J}) and (\ref{dom}) one readily sees that $\Om_J$ is closed
if and only if (\ref{2'}) holds.
In order to determine the \ka nullity ${\cal D}$ we consider the
$J$-anti-invariant 2-forms
$$\phi = e^{\frac{u}{2}}(wdz\wedge dx  + \om\wedge dy),$$
$$J\phi = e^{\frac{u}{2}}(wdz\wedge dy  -\om\wedge dx).$$
They are both of square-norm 2 and we then have
$$d\phi = \tau_{\phi}\wedge \phi; \ \ d (J\phi) = \tau_{J\phi}\wedge
J\phi,$$
where, according to (\ref{na-phi}), the 1-forms $\tau_{\phi},
\tau_{J\phi}$ are given by
\begin{equation}\label{tau-phi}
\tau_{\phi}= -Jb - J\phi(a); \ \  \ \tau_{J\phi}= -Jb + J\phi(a).
\end{equation}
On the other hand, computing $d \phi$ and $d (J\phi)$ directly by
making use of (\ref{dom}) we get
$$\tau_{\phi}= \frac{du}{2} + 2 (\ln w)_y dy; \ \ \tau_{J\phi} =
 \frac{du}{2} + 2 (\ln w)_x dx.$$
We conclude by (\ref{tau-phi}) that
$ J\phi(a)=  (\ln w)_x dx - (\ln w)_y dy.$
But we know from (\ref{na-om}) that $J\phi(a)$ belongs to ${\cal D}$;
the latter  implies the following relations:
\begin{enumerate}
\item[(a)] $(g,J)$ is \ka if and only if $w$ does not depend on $x$
and $y$;
\item[(b)] if $(g,J)$ is not K{\"a}hler, then ${\cal D}= span \{
\frac{\partial}{\partial x},
\frac{\partial}{\partial y} \}$;
\item[(c)] $|\na (\Om_J)|_g^2 = \frac{w_x^2
+ w_y^2}{4e^u w^3}$.
\end{enumerate}
The Ricci form  of the \ka structure $(g,I)$ is given by
$\frac{1}{2}dd^c_{I}u$ (see \cite{LeB}). Here, and in the rest of the
paper, the operator
$d^c_I$ denotes the composition $ I\circ d $, where $d$ is the usual
differential.
Clearly, the Ricci tensor of $g$ is
$J$ invariant if and only if $dd^c_I u$ is a $(1,1)$-form with respect to $J$.
One easily checks that
the latter is equivalent to
$$\big(\frac{u_z}{w}\big)_x= \big(\frac{u_z}{w}\big)_y=0.$$ Thus $u_z
= f w$ for some function $f$ of $z$.
By (\ref{2'}) we get moreover $w=\frac{1}{F+ h},$ where $F$ is a
primitive of $f$, i.e.,
$\frac{d}{dz}F= f$, and $h$ is a function of $x$ and $y$. According to
the relation (a), we know that $h$ is constant if and only if $(g,J)$
is K{\"a}hler. Substituting
into (\ref{1}) we obtain that if  $h$ is not constant, then  $F$ is
constant, or equivalently, $w_z=0, \
u_z=0.$ Thus, if $(g,J)$ is not K{\"a}hler, then  $u$ and $w$ are
functions of $x$ and $y$ and the
equation (\ref{1})  simply means that $w$ is a harmonic function of
$x$ and $y$. The
Ricci tensor is then given by
$${\rm Ric}=(u_{xx} + u_{yy})[dx^2 + dy^2].$$
Therefore, according to Corollary \ref{cor1}, the implication in (b) gives
$(g,J) \in {\cal AK}_3$, while according to Lemma
\ref{lem3}, the equality stated in (c) shows that
$(g,J) \in {\cal AK}_2$ if and only if
$e^u = const.\frac{w_x^2 + w_y^2}{w^3}.$ $\bx$
\begin{cor} \label{cor2}
The inclusions ${\cal K} \subset {\cal AK}_2 \subset {\cal AK}_3$
 are strict in any dimension $2n, n\ge 2$.
\end{cor}
{\it Proof:} Multiplying the examples obtained via Proposition
\ref{prop1} by Riemann surfaces one provides appropriate examples in
any dimension. $\bx$

\section{Classification results}

The proof of Theorem \ref{th1} stated in the introduction
will be a consequence of a slightly more general classification
that we shall prove in Theorem \ref{th2} (see below).
The key idea of the proof is to investigate the properties of
the negative almost complex structure that we define as follows:

\vspace{0.1cm}
\noindent
{\bf Definition.} Let $(M,g,J)$ be a strictly almost-\ka 4-manifold.
On the open set of points where the Nijenhuis tensor of $(g, J)$
does not vanish, let $I$ be the almost complex
structure defined to be equal to $J$ on ${\cal D}$ and to $-J$
on ${\cal D}^{\perp}$.

\vspace{0.1cm} \noindent Clearly, the almost complex structure
$I$ is $g$-orthogonal and yields on the manifold the opposite
orientation than the one given by $J$. We show that curvature
symmetry properties of the almost \ka structure $(g,J, \Om)$ have a
strong effect on the negative almost Hermitian structure
$(g,I,{\bar \Om})$, where ${\bar \Om}$ denotes the fundamental
form of $(g, I)$.

\vspace{0.1cm}

Let us assume that $(M, g, J, \Om)$ is a 4-dimensional, strictly almost \ka
manifold of the class ${\cal AK}_3$. We use the same notations as in
the previous
sections, in particular
for the 1-forms $a$ and $b$ defined by (\ref{na-om}) and (\ref{na-phi})
under the
same convention for the choice of the gauge $\phi$. Our first goal is
to show that the negative almost Hermitian structure $(g, I, {\bar
\Om})$
is almost K{\"a}hler, and then to determine
the 1-forms ${\bar a}, {\bar b}$ corresponding to
the negative gauge
\begin{equation} \label{bargauge}
{\bar \phi} = \phi + \frac{12}{(\kappa - s)} Ja \wedge J\phi(a),
\end{equation}
see (\ref{kappa}).
This is summarized in the following
\begin{Lemma}\label{lem4} Let $(M,g,J,\Om)$ be a strictly almost \ka
4-manifold in the
class ${\cal AK}_3$ and let $I$ be the negative, orthogonal, almost complex
structure defined as above. Then $(g,I, {\bar \Om})$ is an almost \ka
structure compatible
with the reversed orientation of $M$. Moreover,
${\cal D}^{\perp}$ belongs to
the \ka nullity of $(g,I)$ and, with the choice of the negative gauge
as above,
\begin{equation} \label{bar-b}
{\bar b} = 3b + \frac{12\lambda}{(\kappa - s)} \phi(a) .
\end{equation}
\end{Lemma}
{\it Proof:} Defining the 1-forms $m_i, n_i, \ i=1,2$, by
\begin{equation}\label{na-a-0}
\na a = m_1\otimes a + n_1\otimes Ja + m_2\otimes \phi(a) +
n_2\otimes J\phi(a),
\end{equation}
 we use (\ref{na-om}) and
(\ref{na-phi}) to derive the next three equalities:
\begin{eqnarray}\nonumber
\na(Ja) &=& -n_1 \otimes a + m_1 \otimes Ja + (a-n_2)\otimes \phi(a)
\\ \nonumber
& & + (m_2 - Ja) \otimes J\phi(a); \\  \label{na-a-1}
\na (\phi(a)) &=& -m_2 \otimes a + (n_2 - a) \otimes Ja + m_1 \otimes
\phi(a) \\ \nonumber
& & + (b - n_1) \otimes J\phi(a); \\ \nonumber
\na (J\phi(a)) &=& -n_2 \otimes a + (Ja - m_2) \otimes Ja +
(n_1-b) \otimes \phi(a) \\ \nonumber
 & & + m_1 \otimes J\phi(a).
\end{eqnarray}
From (\ref{na-a-0}), (\ref{kappa}) and Lemma
\ref{lem3}-(\ref{lem3-1}) we obtain
\begin{eqnarray}\nonumber
m_1 &=& \frac{1}{|a|^2}g(\na a,a) = \frac{1}{2}d(\ln{(\kappa -s)})\\
\label{m1}
    &=&- \frac{6 \la}{(\kappa -s)} J\phi(a).
\end{eqnarray}
We further use the Ricci relations (\ref{*}) in order
to
determine the 1-forms $n_1, m_2$, and $n_2$.
For that we replace the left-hand sides of the two equalities (\ref{*})
respectively by
$$da = m_1\wedge a + n_1\wedge Ja + m_2\wedge \phi(a) + n_2\wedge
J\phi(a),$$
$$d(Ja) = -n_1 \wedge a + m_1 \wedge Ja + (a-n_2)\wedge \phi(a) +
(m_2 - Ja) \wedge J\phi(a),$$
(see (\ref{na-a-1})), and also take into account that under the ${\cal
AK}_3$ assumption we have
$$R(\phi)= (\frac{s-\kappa}{12} + \lambda)\phi; \ \
R(J\phi)=(\frac{s-\kappa}{12} - \lambda)J\phi,$$
see Lemma \ref{lem1} and (\ref{u(2)})--(\ref{w^+3}).
After comparing the components of both sides, we obtain
\begin{equation}\label{part}
n_1 = -b -\frac{6\la}{(\kappa -s)}\phi(a); \  m_2 = \frac{1}{2}Ja + Jm_0;
\ n_2 =\frac{1}{2}a + m_0,
\end{equation}
where $m_0$ is a 1-form which belongs to ${\cal D}$.

With relations (\ref{na-a-0})--(\ref{part}) in hand, we
can now compute $\na {\bar \Om}$, starting from ${\bar \Om} = \Om -
\frac{12}{(\kappa - s)} a \wedge Ja $ (see (\ref{kappa})), and also
using (\ref{na-om}). We get:
\begin{equation} \label{na-bar-Om}
\na {\bar \Om} = 2m_0 \otimes {\bar \phi} -
2{I}m_0 \otimes {I}{\bar \phi} .
\end{equation}
This proves that $(g, I, {\bar \Om})$ is an almost \ka structure, since
$d{\bar \Om} = 0$ is immediate from (\ref{na-bar-Om}).
The claim about the
\ka nullity of $(g,I)$ follows from  ${\bar a} = 2 m_0 \in {\cal D}$.
Similarly, starting from (\ref{bargauge}) and using (\ref{na-phi}),
(\ref{na-a-0})--(\ref{part}) we obtain
\begin{equation} \label{na-bar-phi}
 \na {\bar \phi} = (3b + \frac{12\lambda}{(\kappa -s)} \phi(a))
\otimes {I}{\bar \phi}
- 2m_0 \otimes {\bar \Om} ,
\end{equation}
and the relation (\ref{bar-b}) follows. $\bx$

\vspace{0.2cm}
As our statements are purely local,
for brevity purposes, we now introduce the following

\vspace{0.1cm}
\noindent
{\bf Definition.} Let $(M,g,J)$ be a strictly almost \ka 4-manifold
in the class ${\cal AK}_3$, and suppose that the Nijenhuis tensor
of $(g,J)$ does not vanish anywhere.
We say that $(M,g,J)$ is a
{\it doubly ${\cal AK}_3$} manifold, if the almost \ka structure $(g,I)$
defined above belongs to the class ${\cal AK}_3$ as well.

\vspace{0.1cm}

\noindent
{\bf Remark 2.} Every non-K{\"a}hler 4-manifold in the class ${\cal
AK}_3$, which is
Einstein, or belongs to class ${\cal AK}_2$ is a doubly ${\cal AK}_3$
manifold. Indeed, this is an immediate consequence of Lemma \ref{lem2}
and Corollary \ref{cor1}.
Note also that all the examples arising
from Proposition 1 are doubly ${\cal AK}_3$ manifolds --- the
negative almost \ka
structure $(g,I)$ is in fact K{\"a}hler for all these examples.

\vspace{0.2cm}

To anticipate, the end result of this section, slightly more general
than Theorem \ref{th1}, will be that every
non-K{\"a}hler, doubly
${\cal AK}_3$ 4-manifold is necessarily given by Proposition \ref{prop1}.
Getting closer to this goal, we now prove
\begin{prop}\label{prop2} Let $(M,g,J)$ be a non-K{\"a}hler,
doubly ${\cal AK}_3$
4-manifold. Then the negative almost \ka structure $(g,I)$ is
K{\"a}hler. Moreover, the Ricci tensor is given by
$$ {\rm Ric} = \frac{s}{2}g^{\cal D}, $$
where $g^{\cal D}$ denotes the restriction
of the metric to the \ka nullity ${\cal D}$ of $(g,J)$. \end{prop}
\noindent
{\it Proof of Proposition \ref{prop2}:}
 For the beginning, we assume only that
$(M, g, J)$ is a strictly almost \ka manifold of the class ${\cal AK}_3$.
We use the Bianchi identity (\ref{lem3-1}), together with
(\ref{lem3-2}) rewritten as
\begin{equation}\label{lem3-2'}
d\lambda = 2 \lambda Jb - \frac{\kappa}{4} J\phi(a) + J\phi({\rm
Ric}_0(a)),
\end{equation}
and the relation (see (\ref{na-a-1})--(\ref{part}))
\begin{equation}\label{na-J-phi-a}
d(J\phi(a)) = - 2b\wedge \phi(a) - m_0\wedge a - Jm_0\wedge Ja .
\end{equation}
Differentiating (\ref{lem3-1}), we get by (\ref{lem3-2'}) and
(\ref{na-J-phi-a}):
\begin{equation} \label{frob1}
0 = 2\lambda (b \wedge \phi(a) - Jb \wedge J\phi(a))
   + \lambda ( m_0 \wedge a + Jm_0 \wedge Ja)
\end{equation}
$$ - J\phi({\rm Ric}_0(a)) \wedge J\phi(a). $$
Taking various components, the relation (\ref{frob1}) can be seen to be
equivalent to:
\begin{equation} \label{nicefrob1}
\lambda m_0 = 2\lambda \phi( b^{{\cal D}^{\perp}}) =
    \frac{1}{2} ({\rm Ric}_0(a))^{{\cal D}} ,
\end{equation}
where the super-scripts ${\cal D}$ and ${\cal D}^{\perp}$ denote the
projections
on those spaces.
Now we shall consider separately the following two cases:

\vspace{0.2cm}
\noindent
{\it Case 1.} $(M,g,J)$ is a doubly ${\cal AK}_3$ manifold which does
not belong to ${\cal AK}_2$. Then by Corollary \ref{cor1} we have
$\la \neq 0$. Since, by assumption, the Ricci tensor is both $J$ and $I$
invariant, it follows that ${\cal D}$ and ${\cal D}^{\perp}$ are eigenspaces
for the traceless Ricci tensor ${\rm Ric}_0$. In other words, we have
\begin{equation}\label{ric0}
{\rm Ric}_0 = \frac{f}{4}[-g^{\cal D} + g^{{{\cal D}}^{\perp}}],
\end{equation}
where $f$ is a smooth function. This implies that $({\rm
Ric}_0(a))^{{\cal D}} = 0$. Since
$\lambda \neq 0$, from (\ref{nicefrob1}) it follows that $m_0 = 0$, {\it i.e.},
$(g, I)$ is K{\"a}hler, see (\ref{na-bar-Om}). Also, from
(\ref{nicefrob1}) it follows that
$ b \in {\cal D}$. Under the doubly ${\cal AK}_3$ assumption, the
Ricci relation (\ref{**}) takes
the form
$$db = a\wedge Ja - \frac{(s+2\kappa)}{12}\Om + \frac{f}{4}{\bar
\Om},$$
or further (see (\ref{kappa}))
\begin{equation}\label{**''}
db =  - \frac{(s+f)}{4}A \wedge JA
           + \frac{(3f-s-2\kappa)}{12} B\wedge JB ,
\end{equation}
where $\{B, JB \}$ is an orthonormal basis for ${\cal D}$ and $\{ A, JA \}$
is an orthonormal basis for ${\cal D}^{\perp}$.
Similarly, the Ricci relation (\ref{**}), written with respect to the
\ka structure $(g,I)$, reads as
\begin{equation}\label{**'''}
d{\bar b} = \frac{(f+s)}{4}A\wedge JA + \frac{(f-s)}{4}B\wedge JB.
\end{equation}
On the other hand, using Lemma \ref{lem3}-(\ref{lem3-1}), the
equality (\ref{bar-b}) can be rewritten as
$${\bar b}  = 3b + d^c_J\ln(\kappa -s),$$
where, we recall, $d^c_J = J \circ d$.
After differentiating we obtain the gauge
independent equality
\begin{equation}\label{db-db}
d{\bar b}=  3db + dd^c_J(\ln(\kappa
-s)).
\end{equation}
For computing $dd^c_J(\ln (\kappa -s))$, we remark first that by
Lemma \ref{lem3}-(\ref{lem3-1}) the vector field dual to
$d^c_J(\ln (\kappa -s))$ belongs to the kernel ${\cal D}$ of the
Nijenhuis tensor of $J$, so that $dd^c_J(\ln (\kappa -s))$ is a
(1,1)-form with respect to $J$. Furthermore, from Lemma
\ref{lem3}-(\ref{lem3-1}) it also follows that $d^c_J\ln(\kappa
-s) = d^c_I (\ln (\kappa -s))$, and then
\begin{equation}\label{explain}
dd^c_J(\ln (\kappa -s))=
 dd^c_I (\ln (\kappa -s)),
\end{equation}
where $d^c_I = I \circ d$ stands for the $d^c$ operator with
respect to $I$. Since $I$ is integrable,  the latter equality
shows that the 2-form $dd^c_J(\ln (\kappa -s))$ it is of type
(1,1) with respect to $I$ as well. Finally, keeping in mind that
$I$ is \ka and $J$ is almost \ka, from (\ref{explain}),
(\ref{lem3-3}) and (\ref{lem3-1}) we compute
\begin{eqnarray}\nonumber
\langle dd^c_J(\ln
(\kappa -s)), {\bar \Om}\rangle &=& \langle dd^c_J(\ln (\kappa -s)), \Om
\rangle  \\ \nonumber
 &=& -\Delta \ln(\kappa -s) \\ \nonumber &=& -\frac{\Delta (\kappa
-s)}{(\kappa -s)}
+\frac{|d(\kappa -s)|^2}{(\kappa -s)^2} \\ \nonumber
 &=& \frac{\kappa -f}{2}.
\end{eqnarray}
Since $dd^c_J(\ln(\kappa -s))$ is a (1,1)-form with respect to both $J$
and $I$, the latter equality shows that
\begin{equation}\label{ddc}
dd^c_J(\ln(\kappa -s))= \frac{(\kappa -f)}{2}
B\wedge JB.
\end{equation}
By (\ref{**''}),  (\ref{**'''}) and (\ref{ddc}), the equality
(\ref{db-db}) finally reduces to $f+s=0$ which, together with
(\ref{ric0}), imply the claimed expression of the Ricci tensor.

\vspace{0.2cm}
\noindent
{\it Case 2.} $(M,g,J)$ is non-K{\"a}hler manifold in the class ${\cal
AK}_2$. Now $\la =0$ by Lemma \ref{lem1},
so the equality (\ref{nicefrob1}) is not useful anymore, as all terms
vanish trivially. However, applying Case 1 to the structure $(g, I)$, we
conclude that it must be itself in the class
${\cal AK}_2$, since otherwise it would follow that $(g, J)$ is K{\"a}hler,
a contradiction.
With the same choices of the gauge as in Lemma \ref{lem4},
we have in this case ${\bar b} = 3b$. This leads to the gauge independent
relation $d{\bar b} = 3db$. Assuming that $(g,I)$ is not K{\"a}hler, we
interchange
the roles of $J$ and $I$ to also get
$db = 3d{\bar b}$, \ {\it i.e.}, $db =0$ holds. But this leads to a
contradiction.
Indeed, according to
Corollary \ref{cor1} we have $f= \kappa$, so from the Ricci relation
(\ref{**''}) we get
$\kappa - s = 0$, {\it i.e.}, $(g,J)$ is K{\"a}hler which contradicts the
assumption.
Thus $(g,I)$ must be K{\"a}hler and
(\ref{**'''}) holds.  It is easily checked that $d{\bar b} = 3db$ is,
in this case,
equivalent to $\kappa + s = 0$.
This and Corollary 1 imply the desired form of the Ricci tensor. $\bx$

\begin{prop}\label{prop4} Let $(M,g,J)$ be a non-K{\"a}hler,
doubly ${\cal AK}_3$ 4-manifold. Then ${\cal D}^{\perp}$ is spanned by
commuting Killing vector fields.
\end{prop}
\noindent
{\it Proof of Proposition \ref{prop4}:}
For any  smooth functions $p$ and
$q$ we consider the vector field $X_{p,q}$ in ${\cal D}^{\perp}$, the dual to
the 1-form $pa + qJa$. The condition that $X_{p,q}$ is Killing
is equivalent to $\na (pa + qJa)$ being a section of $\La^2M$. To
write explicitly the equation on $p$ and $q$ that arise from the
latter condition we need the covariant derivative of $a$
and $Ja$. But we know
already
by Proposition \ref{prop2} that $(g,I)$ is K{\"a}hler, {\it i.e.}, the 1-form
$m_0$ defined in (\ref{part})  vanishes (see (\ref{na-bar-Om})).
We thus have by (\ref{na-a-0})--(\ref{part})
\begin{eqnarray}\label{na-a}
\na a &=& -\frac{6\la}{(\kappa -s)}J\phi(a)\otimes a -
\frac{6\la}{(\kappa -s)} \phi(a)\otimes Ja \\ \nonumber
        & & - b\otimes Ja + \frac{1}{2}Ja\otimes \phi(a) +
\frac{1}{2}a\otimes J\phi(a); \\ \label{na-J-a}
\na (Ja) & = & \frac{6\la}{(\kappa - s)} \phi(a)\otimes a -
 \frac{6\la}{(\kappa -s)} J\phi(a)\otimes Ja \\ \nonumber
         & &  + b\otimes a + \frac{1}{2} a\otimes\phi(a) + \frac{1}{2}
Ja\otimes J\phi(a).
\end{eqnarray}
Using (\ref{na-a}) and (\ref{na-J-a}) the
condition that $\na(pa + qJa)$ belongs to $\La^2M$  can be
rewritten as
\begin{equation}\label{system}
\left.
\begin{array}{c@{ \ = \ }c}
dp & -qb  - \frac{p}{2}(1-\frac{12\la}{(\kappa - s)})J\phi(a) -
\frac{q}{2}(1 + \frac{12\la}{(\kappa -s)})\phi(a)  \ + rJa; \\
dq & pb  - \frac{p}{2}(1-\frac{12\la}{(\kappa - s)}) \phi(a)  \  +
\frac{q}{2}(1 + \frac{12\la}{(\kappa -s)})J\phi(a) - ra,
\end{array}
\right\}
\end{equation}
where $r$ is a smooth function. Since we are looking for
commuting Killing fields, we have $r\equiv 0$, and we thus obtain
a Frobenius type system. To show that (\ref{system}) has solution
in a neighborhood of a point $x\in M$ for any given values
$(p(x), q(x))$, we apply the Frobenius theorem. Accordingly, we
have to check
\begin{equation}\label{edno}
d\big( \ \ 2qb + p(1-\frac{12\la}{(\kappa - s)})J\phi(a) +
q(1 + \frac{12\la}{(\kappa -s)})\phi(a)\big)= 0;
\end{equation}
\begin{equation}\label{dve}
d\big(-2pb + p(1-\frac{12\la}{(\kappa - s)})\phi(a) - q(1 +
\frac{12\la}{(\kappa -s)})J\phi(a)\big)=0.
\end{equation}
For that we further specify the relations (\ref{na-a-1}) and
(\ref{**''}), taking into account that  $m_0=0$ and
$f=-s$ (see Proposition \ref{prop2}). We thus get:
$$d(J\phi(a)) = -2 Jb \wedge J\phi(a),$$
$$d(\phi(a)) =  2b\wedge J\phi(a)+ 2\la B\wedge JB,$$
$$db = -\frac{(2s+ \kappa)}{6} B\wedge JB,$$
where $B=\frac{1}{|a|}\phi(a)$ and $JB=\frac{1}{|a|}J\phi(a)$ is an
orthonormal frame of ${\cal D}$. By Lemma \ref{lem3} and  (\ref{ddc})
we also have
$$ d\ln(\kappa -s) = -\frac{12\la}{(\kappa - s)}J\phi(a); \ \
d^c_J\ln(\kappa -s)= \frac{12\la}{(\kappa -s)}\phi(a),$$
$$
dd^c_J(\ln(\kappa -s))= \frac{(\kappa +s)}{2}B\wedge JB.$$
Using the above equalities, together with (\ref{system}) and
(\ref{kappa}), it is
now straightforward to check (\ref{edno})
and (\ref{dve}). $\bx$

\vspace{0.2cm}
\noindent
{\bf Remark 3.} The miraculous
cancellation that appears by checking the equalities (\ref{edno}) and
(\ref{dve}) can be explained by simply observing that if the cancellation
hadn't occurred we would then derive an integrability condition depending
on $\lambda $ and $\kappa -s$. But these take arbitrary values for the
examples provided by Proposition 1. We thus conclude that the integrability
conditions (\ref{edno}) and (\ref{dve}) must be satisfied.

\begin{theo}\label{th2} Any 4-dimensional non-K{\"a}hler,  doubly ${\cal
AK}_3$ metric
 is locally isometric to one of the metrics described by Proposition
{\rm \ref{prop1}(i)} {\rm (}or equivalently, by {\rm (\ref{ak3})}{\rm )}.
\end{theo}
\noindent {\it Proof of Theorem 2:} Let $(M,g,J)$ be a
non-K{\"a}hler,  doubly ${\cal AK}_3$ 4-manifold. By Proposition
\ref{prop2}, there exists a \ka structure $I$, which yields the
opposite orientation of $M$. Moreover, we know by Proposition
\ref{prop4} that in a neighborhood of any point there exists a
Killing vector field $X \in {\cal D}^{\perp}$, determined by a
solution of the system (\ref{system}). It is not difficult to
check that $X$ preserves $I$. Indeed, we have to verify
$${\cal L}_X {\bar \Om} = d(I(pa + q Ja)) = d(qa -pJa)=0.$$
The latter equality is a consequence of (\ref{system})
and the Ricci identities (\ref{*}) (If the manifold is not Ricci flat,
the invariance of $I$ also follows from the fact that
$I$ is determined up to sign  by the two eigenspaces of ${\rm Ric}$).
According to
\cite{LeB}, the metric
$g$ has the form (\ref{g}), where the functions $w$ and $u$ satisfy
(\ref{1}) and
$X=\frac{\partial}{\partial t}$. From Proposition \ref{prop2}  we also know
that ${\rm Ric}(X)=0$. But the Ricci form of the \ka structure $(g,I)$
is given by
$\frac{1}{2}dd^c_{I}u$ (see \cite{LeB}); we thus obtain $w=const.u_z$ and
then
$${\rm Ric}=(u_{xx} + u_{yy} + (e^u)_{zz})[dx^2 + dy^2].$$
The above equality shows that either $g$ is Ricci flat (then $g$ is
given by Tod's ansatz, see \cite{Arm1}),
or else, according to Proposition \ref{prop2},
the \ka nullity ${\cal D}$ of $(g,J)$ is spanned by the (Riemannian) dual
fields
$dx$ and  $dy$. The latter means that the \ka \  form $\Om$ of $(g,J)$
is given by (\ref{J}), and the result follows by Proposition \ref{prop1}
and Remark 1.
$\bx$

\vspace{0.2cm}

\noindent
Theorem 1 is now just a particular case.

\vspace{0.2cm}

\noindent
{\bf Proof of Theorem 1.} By Remark 2 we know that every strictly
almost-K\"ahler 4-manifold $(M,g,J,\Om)$ satisfying ($G_2$) is
doubly ${\cal AK}_3$;
it follows by Theorem 2 and Proposition 1 that $(M,g,J,\Om)$ arises from
Proposition
\ref{prop1}(ii). According to Remark 1(b) the metric $g$ is
locally isometric to
(\ref{canonic}) which, in turn, is isometric to Kowalski's metric, doing
the change
of variables (\ref{chvar}). $\bx$

\vspace{0.1cm}
\noindent
{\bf Remark 4.} Avoiding the use of the change of variables (\ref{chvar}),
one could have
completed the proof of Theorem 1 as follows: as above one shows that
any strictly almost-K\"ahler 4-manifold $(M,g,J,\Om)$ satisfying ($G_2$) is
locally isometric to (\ref{canonic}). On the other hand,
A.Gray \cite{gray} showed that any Riemannian  3-symmetric
space has a canonical almost-Hermitian structure, which in 4-dimensions, is
necessarily almost-K{\"a}hler (K{\"a}hler iff the manifold is
symmetric) and satisfies the condition  ($G_2$).
It thus follows that the proper 3-symmetric metric of Kowalski \cite{kow}
is isometric to (\ref{canonic}) as well.  In particular, this
provides a differential geometric proof of existence and uniqueness
of proper
3-symmetric 4-dimensional manifolds, result proved by Kowalski
using Lie algebra techniques \cite{kow}.

\begin{cor}\label{cor3} {\rm { (\cite{ADK})}} Every compact almost
\ka 4-manifold satisfying the second curvature condition of Gray is
K{\"a}hler.
\end{cor}
{\it Proof of Corollary \ref{cor3}:} Suppose for contradiction that
$(M,g,J)$ is a compact,
non-K{\"a}hler, almost \ka 4-manifold in the class ${\cal AK}_2$.

According to Corollary \ref{cor1},
the distributions ${\cal D}$ and ${\cal D}^{\perp}$ are globally
defined on $M$, and by Proposition 2 they give rise to a
negative \ka structure $(g,I)$. We
know by Theorem 1 that $(g,J,I)$ locally arise from Proposition 1. Then
the whole curvature of $g$ is completely determined
by the (negative constant) scalar curvature $s$, cf. Remark 1. More precisely,
the conformal
curvature $\kappa$ is given by $\kappa= -s$ (Corollary \ref{cor1} and
Proposition \ref{prop2}).  Since
$(g,I)$ is K{\"a}hler, we also have $|W^-|^2 = \frac{s^2}{24}$, see {\it
e.g.} \cite{Ga}. As
$(g,J)$ is in the class ${\cal AK}_2$, the self-dual Weyl tensor satisfies
$W^+_2=0$, $W^+_3=0$ and then $|W^+|^2 = \frac{\kappa^2}{24}$ (see
(\ref{w^+1})); by $\kappa = -s$ we conclude $|W^+|^2 = |W^-|^2
=\frac{s^2}{24}$.
We then get by the Chern-Weil formula
$$\sigma(M)=\frac{1}{12\pi^2}\int_M
|W^+|^2 - |W^-|^2 dV_g$$
that the signature $\sigma(M)$ vanishes.
Similarly, the Euler characteristic $e(M)$ is given
by
$$e(M) = \frac{1}{8\pi^2}\int_M |W^+|^2 + |W^-|^2 + \frac{s^2}{24} -
\frac{1}{2}|{\rm Ric}_0|^2 dV_g.$$
But we know that the Ricci tensor of $g$ has
eigenvalues $(0,0,\frac{s}{2},\frac{s}{2})$ (Proposition \ref{prop2}) and
then $|{\rm Ric}_0|^2=\frac{s^2}{4}$; we thus readily see that $e(M)=0$.
Furthermore, since $({\bar M},g,I)$ is a K{\"a}hler surface of
(constant) negative
scalar curvature, we have $H^0({\bar M}, { K}^{\otimes - m})=0$, where
$K$ denotes the canonical bundle of $({\bar M},I)$. The
conditions
$\sigma({\bar M})=-\sigma(M)=0$, $e({\bar M})=e(M)=0$  then imply that
the Kodaira dimension of $({\bar M}, I)$ is necessarily equal to 1,
cf. {\it e.g.}
\cite{BPV}. Thus $({\bar M},I)$
is a minimal
properly elliptic surface with vanishing Euler
characteristic. Using
an argument from \cite{ADK}, we conclude that, up to a finite cover, $({\bar
M},I)$ admits a
non-vanishing holomorphic vector field $X$. Now the well known
Bochner formula for holomorphic fields and the fact that the
Ricci tensor of $({\bar M},g,I)$ is semi-negative whose kernel is the
distribution ${\cal D}^\perp$ (Proposition \ref{prop2}) imply
that $X$ is
parallel and belongs to  ${\cal D}^{\perp}$. Then ${\cal
D}^{\perp}$ (hence
also ${\cal D}$) is parallel. Since $(g,I)$ is a K{\"a}hler structure,
$I$ is parallel, and consequently,  the almost
complex structure
$J$ must be parallel as well, {\it i.e.}, $(g,J)$ is K{\"a}hler,
which contradicts our assumption. $\bx$

\vspace{0.1cm}
\noindent
{\bf Remark 5.} For obtaining a contradiction
in the proof of Corollary \ref{cor3}
one can  alternatively argue as follows: We
know by Theorem 1 that $(g,J,I)$ locally arise from Proposition 1.
The metric $g$ is therefore locally homogeneous and the complex structure $I$
is invariant  as being determined by the eigenspaces of the Ricci
tensor. It thus follows that $(M,g,I)$ is a {\it compact} locally homogeneous
K\"ahler surface; it is well known that any such  surface is locally
(Hermitian) symmetric (cf. {\it e.g.} \cite{wall}), while  the metric
$g$ given by Proposition 1(ii) is not.

\vspace{0.1cm}
\noindent
{\bf Remark 6.} Using the method of ``nuts and bolts'' \cite{GH},
C.LeBrun \cite{LeB1} successfully ``compactified'' certain K{\"a}hler
metrics arising from (\ref{g}) and obtained explicit examples
of compact scalar-flat K{\"a}hler surfaces admitting a circle
action. The idea is the following: Starting from an open (incomplete)
manifold $M_0$
where the metric $g$ has the form (\ref{g}), one adds points and (real)
surfaces in order to obtain a larger, complete manifold $M$, such that $M_0$ is
a dense open subset of $M$, and the circle action on $M_0$ generated
by the Killing vector field $X=\frac{\partial}{\partial t}$ extends to
$M$; the added points and surfaces become the fixed point of this
action.

It is thus natural to wonder if similar ``compactification'' exists
for the metrics given by Proposition \ref{prop1}, providing
compact
examples of non-K{\"a}hler, almost \ka
4-manifolds in the class ${\cal AK}_3$. (The interest in such compact examples
is motivated by some variational problems on compact symplectic
manifolds \cite{Bl,BI}). Corollary \ref{cor3} shows that this is
impossible if we
insist that (\ref{u}) is satisfied. Unfortunately, even in the case when
$(\ref{u})$ does not hold, the variable reduction we have for the
functions $u$ and $w$ does not permit us to obtain compact
examples directly following LeBrun's approach. Indeed, if $(M,g,J)$ was
a compactification of $(M_0,J,g)$ with
extended circle action  generated by a Killing vector field
$X= \frac{\partial}{\partial t}$, then by
Propositions \ref{prop2} and \ref{prop4}, we would have
${\rm Ric}(X,X)=0$ on $M_0$, hence
also, on $M$ as $M_0$ is a dense subset; by the Bochner formula
$X$ would then be parallel. In particular, the $g$-norm of $X$ would be
constant, hence also, the smooth function $w=\frac{1}{g(X,X)}$. Therefore,
$(g,J)$ would be K{\"a}hler by Proposition \ref{prop1}, a contradiction.

\vspace{0.2cm} \noindent As a final note, it is tempting to
conjecture that the local classification obtained in Theorem 2
could be further extended to the general case of strictly ${\cal
AK}_3$ 4-manifolds (in other words, we believe that the doubly
${\cal AK}_3$ assumption in the Theorem 2 could be removed). For
this goal a further analysis of the higher jets of $J$ would be
needed, with computations becoming more involved, but it is
possible that some nice cancellations might still take place.

\vspace{0.2cm}
\noindent  {\bf Acknowledgements}: The first author thanks the
Mathematical Institute of Oxford for hospitality during the preparation
of an important part of this paper.
The authors are
grateful to R. Bryant, G. Gibbons, C. LeBrun, S. Salamon and P. Tod for
their interest and  some stimulating discussions. We would also like to
express our thanks to O. Mu\u{s}karov whose comments essentially improved
the presentation of the results in Section 4, to D. Blair for his friendly
assistance in reading the manuscript and suggesting several improvements,
and to A. Moroianu for bringing to our attention the unpublished work
\cite{Br}.

\end{document}